\documentclass[12pt]{amsart}
\usepackage{amsmath, amssymb, latexsym, amsthm, mathrsfs, geometry}
\usepackage[all]{xy}

\newcommand{\ed}{

\subsection*{Acknowledgments}

We thank Salvador Herna\'ndez for inviting the second named author to give a lecture on the
topic in the conference \emph{Functional Analysis in Valencia 2010},
dedicated to the 80'th birthday of Manuel Valdivia.
We also thank Lyubomyr Zdomskyy for reading the paper and making useful comments.

This work is an extension of a part of the first named author's M.Sc.\ thesis
at the Weizmann Institute of Science, supervised by Gady Kozma and the second named author.
We thank Gady Kozma for useful discussions, and the Weizmann Institute of Science for the stimulating atmosphere.

\end{document}}

      \newenvironment{changemargin}[2]{\begin{list}{}{
         \setlength{\topsep}{0pt}\setlength{\leftmargin}{0pt}
         \setlength{\rightmargin}{0pt}
         \setlength{\listparindent}{\parindent}
         \setlength{\itemindent}{\parindent}
         \setlength{\parsep}{0pt plus 1pt}
         \addtolength{\leftmargin}{#1}\addtolength{\rightmargin}{#2}
         }\item }{\end{list}}
\newcommand{\nc}{\newcommand}
\nc{\Om}{\Omega}
\nc{\Ga}{\Gamma}
\nc{\ga}{\gamma}
\nc{\vphi}{\varphi}
\renewcommand{\liminf}{\op{liminf}}
\nc{\one}{\mathbf{1}}
\nc{\CX}{\op{C}(X)}
\nc{\CY}{\op{C}(Y)}
\nc{\BX}{\op{B}(X)}
\nc{\BaireOneX}{\op{Baire}_1(X)}
\nc{\LI}{\op{LI}}
\nc{\lims}[1]{\cl{#1}^\text{limits}}
\nc{\partlims}{\op{partlims}}
\newcommand{\Pa}[8]{\bibitem{#1} {#2}, \emph{#3}, {#4} \textbf{#5} ({#6}), {#7}--{#8}.}

\newcommand{\Bc}[9]{\bibitem{#1} {#2}, \emph{#3}, in: \textbf{#4} (#5), #6 #7, #8--#9.}

\newcommand{\Setting}[7]{\xymatrix@R=4pt@C=7pt{#1\ar@{-}[r]&#2\ar@{-}[r]&#3\\&#4\ar@{-}[u]\\
#5\ar@{-}[uu]\ar@{-}[r] & #6\ar@{-}[u]\ar@{-}[r] & #7\ar@{-}[uu]}}

\newcommand{\bq}{\begin{quote}}
\newcommand{\eq}{\end{quote}}
\newcommand{\cl}[1]{\overline{#1}}
\newcommand{\CH}{the Continuum Hypothesis}
\newcommand{\MA}{Martin's Axiom}

\newcommand{\inv}{^{-1}}

\newcommand{\sr}[2]{{\txt{$#1$\\#2}}}

\newcommand{\N}{\mathbb{N}}
\newcommand{\NN}{{\N^{\N}}}

\newcommand{\roth}{{[\N]^{\oo}}}
\newcommand{\Fin}{{[\N]^{<\oo}}}

\newcommand{\sseq}[1]{\{#1 : n\in\N\}}

\newcommand{\op}{\operatorname}

\newcommand{\scrA}{\mathscr{A}}
\newcommand{\scrB}{\mathscr{B}}

\newcommand{\cB}{\mathcal{B}}

\newcommand{\cA}{\mathcal{A}}

\newcommand{\cF}{\mathcal{F}}

\newcommand{\rmO}{\mathrm{O}}
\newcommand{\rmL}{\mathrm{L}}

\newcommand{\R}{\mathbb{R}}
\newcommand{\cU}{\mathcal{U}}
\newcommand{\Union}{\bigcup}
\newcommand{\cV}{\mathcal{V}}
\newcommand{\cW}{\mathcal{W}}

\newcommand{\Impl}{\Rightarrow}
\long\def\forget#1\forgotten{}
\newcommand{\ft}{\mathfrak{t}}
\newcommand{\fb}{\mathfrak{b}}

\newcommand{\oo}{\infty}

\newcommand{\fp}{\mathfrak{p}}

\newcommand{\comp}{^{\text{\tt c}}}
\newcommand{\nin}{\notin}

\newcommand{\sbst}{\subseteq}

\newcommand{\sm}{\setminus}

\newcommand{\dom}{\op{dom}}

\newtheorem{thm}{Theorem}[section]
\newcommand{\bthm}{\begin{thm}} \newcommand{\ethm}{\end{thm}}
\newtheorem{prop}[thm]{Proposition}
\newcommand{\bprp}{\begin{prop}} \newcommand{\eprp}{\end{prop}}
\newtheorem{fact}[thm]{Fact}
\newcommand{\bfct}{\begin{fact}} \newcommand{\efct}{\end{fact}}
\newtheorem{prob}[thm]{Problem}
\newcommand{\bprb}{\begin{prob}} \newcommand{\eprb}{\end{prob}}
\newtheorem{lem}[thm]{Lemma}
\newcommand{\blem}{\begin{lem}} \newcommand{\elem}{\end{lem}}
\newtheorem{claim}[thm]{Claim}
\newcommand{\bclm}{\begin{claim}} \newcommand{\eclm}{\end{claim}}
\newtheorem{cor}[thm]{Corollary}
\newcommand{\bcor}{\begin{cor}} \newcommand{\ecor}{\end{cor}}
\newtheorem{conj}[thm]{Conjecture}
\newcommand{\bcnj}{\begin{conj}} \newcommand{\ecnj}{\end{conj}}
\theoremstyle{definition}
\newtheorem{defn}[thm]{Definition}
\newcommand{\bdfn}{\begin{defn}} \newcommand{\edfn}{\end{defn}}
\newtheorem{cnv}[thm]{Convention}
\newcommand{\bcnv}{\begin{cnv}} \newcommand{\ecnv}{\end{cnv}}
\theoremstyle{remark}
\newtheorem{rem}[thm]{Remark}
\newcommand{\brem}{\begin{rem}} \newcommand{\erem}{\end{rem}}
\newtheorem{exam}[thm]{Example}
\newcommand{\bexm}{\begin{exam}} \newcommand{\eexm}{\end{exam}}
\newcommand{\bpf}{\begin{proof}} \newcommand{\epf}{\end{proof}}
\newcommand{\be}{\begin{enumerate}}
\newcommand{\ee}{\end{enumerate}}
\newcommand{\bi}{\begin{itemize}}
\newcommand{\itm}{\item}
\newcommand{\ei}{\end{itemize}}

\forget
\forgotten

\newcommand{\sone}{\mathsf{S}_1}
\newcommand{\sfin}{\mathsf{S}_\mathrm{fin}}
\newcommand{\ufin}{\mathsf{U}_\mathrm{fin}}

\title[The Gerlits--Nagy Problem]{Pointwise convergence of partial functions:\\
The Gerlits--Nagy Problem}
\author{Tal Orenshtein}
\address[Orenshtein]{Department of Mathematics,
Weizmann Institute of Science, Rehovot 76100, Israel}
\email{tal.orenshtein@weizmann.ac.il}
\urladdr{http://www.wisdom.weizmann.ac.il/\~{}talo}
\author{Boaz Tsaban}
\address[Tsaban]{Department of Mathematics, Bar-Ilan University, Ramat Gan 52900, Israel}
\email{tsaban@math.biu.ac.il}
\urladdr{http://www.cs.biu.ac.il/\~{}tsaban}

\begin{document}

\begin{abstract}
For a set $X\sbst\R$, let $B(X)\sbst\R^X$ denote the space of Borel real-valued functions on $X$,
with the topology inherited from the Tychonoff product $\R^X$.
Assume that for each countable $A\sbst B(X)$, each $f$ in the closure of
$A$ is in the closure of $A$ under pointwise limits of sequences of partial functions.
We show that in this case, $B(X)$ is countably Fr\'echet--Urysohn, that is,
each point in the closure of a countable set is a limit of a sequence of elements of that set.
This solves a problem of Arnold Miller. The continuous version of this problem
is equivalent to a notorious open problem of Gerlits and Nagy.
Answering a question of Salvador Herna\'ndez, we show that
the same result holds for the space of all Baire class $1$ functions on $X$.

We conjecture that, in the general context, the answer to the
continuous version of this problem is negative, but we identify a nontrivial context
where the problem has a positive solution.

The proofs establish new local-to-global correspondences,
and use methods of infinite-combinatorial topology,
including a new fusion result of Francis Jordan.
\end{abstract}

\maketitle

\section{Introduction and basic results}

Let $X\sbst\R$. $\CX$ is the family of all continuous real-valued functions
on $X$. We consider $\CX$ with the topology inherited from the Tychonoff product
$\R^X$.  A basis of the topology is given by the sets
$$[f; x_1,\dots,x_k; \epsilon]:=\{g\in \CX : (\forall i=1,\dots,k)\ |g(x_i)-f(x_i)|<\epsilon\},$$
where $f\in \CX$, $k\in\N$, $x_1,\dots,x_k\in X$, and $\epsilon$ is a positive real number.
This is the \emph{topology of pointwise convergence}, where a sequence (more generally, a net) $f_n$
converges to $f$ if and only if for each $x\in X$, the sequence of real numbers $f_n(x)$ converges to $f(x)$.

By definition, the (topological) closure $\cl{A}$ of a set $A\sbst\CX$ is the set of all $f\in\CX$ such that,
for all $k\in\N$, $x_1,\dots,x_k\in X$, and positive $\epsilon$, there is an element $g\in A$ such that
$|g(x_i)-f(x_i)|<\epsilon$ for $i=1,\dots,k$. (Equivalently,
there is a net in $A$ converging pointwise to $f$.)
$\CX$ is metrizable only when $X$ is countable, and thus it makes sense to ask,
when $X$ is not countable, when do limits of sequences determine the closure of sets.

For a topological space $Y$ and $A\sbst Y$, the \emph{closure of $A$ under limits of sequences}
is the smallest set $C\sbst Y$ containing $A$, such that
for each convergent (in $Y$) sequence of elements of $C$, the limit of this sequence is
also in $C$.
The closure of $A$ under limits of sequences is contained in the topological closure
$\cl{A}$ of $A$ in $Y$.

Gerlits \cite{Gerlits83}, and independently Pytkeev \cite{Pytkeev82},
proved that if limits determine the closure in $\CX$, then indeed
it suffices to take limits once.

\bthm[{Gerlits, Pytkeev}]\label{GerPyt}
Let $X$ be a Tychonoff space.
Assume that, for each $A\sbst \CX$, each $f\in\cl{A}$ (closure in $\CX$)
belongs to the closure (in $\CX$) of $A$ under limits of sequences.
Then, for each $A\sbst \CX$, each $f\in\cl{A}$ is a limit of a sequence
of elements of $A$.
\ethm

The properties of $\CX$ in the premise and in the conclusion of Theorem \ref{GerPyt}
are often named \emph{sequential} and \emph{Fr\'echet--Urysohn},
respectively.

\medskip

Consider now \emph{partial} functions $f:X\to\R$, that is,
functions whose domain is a (not necessarily proper) subset of $X$.

\bdfn
Let $f_1,f_2,\dots:X\to\R$ be partial functions.
The \emph{partial limit} function $f=\lim_nf_n$ is the partial real-valued function on $X$,
with $\dom(f)$ being the set of all $x$ such that $f_n(x)$ is eventually
defined and converges, defined by $f(x)=\lim_nf_n(x)$ for each $x\in\dom(f)$.
\edfn

Thus, for $f_1,f_2,\dots\in \CX$, the ordinary limit $\lim_nf_n$ exists in $\CX$ if and only if
the domain of the partial limit function $f=\lim_nf_n$ is $X$, and $f$ is continuous.
The partial limit of a sequence of partial functions always exist, though it may be the empty function.

\bdfn
For a set $A$ of partial functions $f:X\to\R$, the
\emph{closure of $A$ under partial limits of sequences},
$\partlims(A)$, is the smallest set $C$ of partial functions $f:X\to\R$, such that
$A\sbst C$ and for each sequence in $C$, the partial limit of this sequence is
also in $C$.
\edfn

Thus, the closure, in $\CX$, of a set $A\sbst \CX$ under limits of sequences is a subset of $\CX\cap\partlims(A)$.

\blem\label{plimsbcl}
For each $A\sbst \CX$, $\CX\cap\partlims(A)$ is contained in $\cl{A}$, the closure of $A$ in $\CX$.
\elem
\bpf
The definition of basic open sets in $\CX$ (or $\R^X$)
may be extended to partial functions, by letting $[f; x_1,\dots,x_k;\epsilon]$
be the set of all partial $g:X\to\R$ such that
$x_1,\dots,x_k\in\dom(g)$ and $|g(x_i)-f(x_i)|<\epsilon$, for all $i=1,\dots,k$.

Assume that $f\nin\cl{A}$. Take $x_1,\dots,x_k\in X$ and $\epsilon>0$, such that
$A\cap [f; x_1,\dots,x_k;\epsilon]=\emptyset$.
Then $A\sbst [f; x_1,\dots,x_k;\epsilon]\comp$, and $[f; x_1,\dots,x_k;\epsilon]\comp$
is closed under limits of partial functions: Assume and $g=\lim_ng_n\in [f; x_1,\dots,x_k;\epsilon]$.
Then $x_1,\dots,x_k\in\dom(g)$, and $|g(x_i)-f(x_i)|<\epsilon$,
and therefore the same holds for $g_n$, for all but finitely many $n$.
In particular, it cannot be the case that $g_1,g_2,\dots\in[f; x_1,\dots,x_k;\epsilon]\comp$.

It follows that $f$ is not in the closure of $A$ under partial limits of sequences.
\epf


In 1982, Gerlits and Nagy published their seminal paper \cite{GN}. 
This paper has generated over 200 subsequent papers and a rich theory.
Among the problems posed in \cite{GN}, only one remains open. 
On its surface, the \emph{Gerlits--Nagy Problem} is a combinatorial one, and we defer its
combinatorial formulation to Section \ref{OnGNP}, 
where we prove that the Gerlits--Nagy Problem is equivalent to the following fundamental problem, 
dealing with pointwise convergence of real-valued functions.



\bprb[Gerlits--Nagy \cite{GN}]\label{GNP}
Assume that, for each $A\sbst \CX$, each $f\in\cl{A}$ belongs to the
closure of $A$ under partial limits of sequences.
Does it follow that, for each $A\sbst \CX$, each $f\in\cl{A}$ is a limit of a sequence
of elements of $A$?
\eprb

In the \emph{Second Workshop on Coverings, Selections, and Games in Topology}
(Lecce, Italy, 2005), Arnold Miller delivered a plenary lecture, where
he posed the variant of the Gerlits--Nagy Problem, dealing with \emph{Borel} rather than
continuous functions \cite{MilProb}.

Let $\BX\sbst\R^X$ be the family of all Borel real-valued functions on $X$.
One may consider the questions discussed above also for $\BX$, with
the following reservation: Here, one must restrict attention to \emph{countable} $A\sbst \BX$,
as we now show.

Each of the properties mentioned in the above discussion implies
that $\CX$ is \emph{countably tight}, that is,
each point in the closure of a set is in the closure of a countable
subset of that set.
The standard proof would be by transfinite induction on the countable ordinals, but we
adopt here an argument given in \cite{ChMaTaAng}.

\bprp
Let $X$ be a topological space.
Assume that, for each $A\sbst \CX$, each $f\in\cl{A}$ belongs to the
closure of $A$ under partial limits of sequences.
Then $\CX$ is countably tight.
\eprp
\bpf
Let $A\sbst\CX$. By Lemma \ref{plimsbcl}, $\partlims(A)\cap\CX\sbst\cl{A}$.
Thus, it suffices to show that for each $f\in\cl{A}$, there is
a countable $D\sbst A$ such that $f\in\partlims(D)$.

Let $B=\Union\{\partlims(D) : D\sbst A\mbox{ is countable}\}$.
Then $B$ is closed under partial limits of sequences:
Let $f_1,f_2,\dots\in B$. Then there are countable $D_1,D_2,\dots\sbst A$,
such that $f_n\in \partlims(D_n)$ for all $n$. Let $D=\Union_n D_n$. Then
$f_1,f_2,\dots\in\partlims(D)$, and therefore $\lim_nf_n\in\partlims(D)\sbst B$.

Thus, $\partlims(A)\sbst B$, as required.
\epf

By a classical result of Arhangel'ski\u{\i}, $\CX$ is countably tight for all $X\sbst\R$ (indeed,
for all topological spaces $X$ such that all finite powers of $X$ are Lindel\"of).
However, $\BX$ is not countably tight, unless $X$ is countable (in which case, $\R^X$, and thus $\BX$, is metrizable).

We denote by $\one$ the constant function identically equal to $1$ on $X$.

\bprp
Let $X$ be an uncountable space, where each singleton is Borel.
Then $\BX$ is not countably tight.
\eprp
\bpf
Take $A=\{\chi_F : F\sbst X\mbox{ finite}\}\sbst \BX$, where $\chi_F$ denotes the
characteristic function of $F$. Then the constant function $\one$ is in $\cl{A}$.
Let $D=\sseq{\chi_{F_n}}\sbst X$. Take $a\in X\sm \Union_nF_n$. Then $\chi_{F_n}(a)=0$ for
all $n$, and thus $\one\nin\cl{D}$.
\epf

\bprb[Miller 2005 \cite{MilProb}]\label{MP}
Assume that, for each countable $A\sbst \BX$, each $f\in\cl{A}$ belongs to the
closure of $A$ under partial limits of sequences.
Does it follow that, for each countable $A\sbst \BX$, each $f\in\cl{A}$ is a limit of a sequence
of elements of $A$?
\eprb

Our main result (Section \ref{MilSol}) is a solution, in the affirmative, of Miller's problem.
At the end of the second named author's talk in the conference \emph{Functional Analysis in Valencia 2010},
Salvador Herna\'ndez asked what is the solution to Miller's Problem when considering
Baire class $1$ functions (i.e., functions which are pointwise limits of sequences of 
continuous functions).
We solve Herna\'ndez's problem in Section \ref{HernSol}.
Finally, we establish several results concerning the original Gerlits--Nagy Problem, and pose some
related problems.

\section{Borel functions (Miller's Problem)}\label{MilSol}

We solve Miller's Problem \ref{MP} in the affirmative. Indeed, we do so
not only for sets $X\sbst\R$, but for all topological spaces $X$.

\bthm\label{main}
Let $X$ be a topological space.
Assume that, for each countable $A\sbst \BX$, each $f\in\cl{A}$ belongs to the
closure of $A$ under partial limits of sequences.
Then for each countable $A\sbst\BX$, each $f\in\cl{A}$ is a limit of a sequence
of elements of $A$.
\ethm

The proof is divided naturally into four steps. For brevity, we make the following convention,
that will hold throughout the paper.

\bcnv
Let $X$ be a topological space. We say that $\cU$ is a \emph{cover}
of $X$ if $X=\Union\cU$, but $X\nin\cU$.
By \emph{Borel cover} of $X$ we always mean a \emph{countable} family $\cU$ of Borel subsets of $X$,
such that the union of all members of $\cU$ is $X$.
\ecnv

\subsubsection*{Step 1: Local to global}
We deduce from the given local property of $\BX$, a global property of $X$.

\bdfn[Gerlits--Nagy \cite{GN}]
A cover $\cU$ of $X$ is an \emph{$\omega$-cover} of $X$ if
each finite $F\sbst X$ is contained in a member of $\cU$.

For sets $B_1,B_2,\dots$, let
$$\liminf_n B_n=\Union_m\bigcap_{n\ge m}B_n,$$
that is, the set of all $x$ which belong to $B_n$ for all but
finitely many $n$.
Let $\LI(\cU)$ be the closure of $\cU$ under the operator $\liminf$.
\edfn
A basic property of $\liminf_n B_n$ is that it
does not depend on the first few sets $B_n$.

\blem\label{l2g}
Let $X$ be a topological space.
Assume that, for each countable $A\sbst \BX$, each $f\in\cl{A}$ belongs to $\partlims(A)$.
Then for each Borel $\omega$-cover $\cU$ of $X$, $X\in\LI(\cU)$.
\elem
\bpf
Let $\cU$ be a Borel $\omega$-cover of $X$.
Take $A=\{\chi_U : U\in\cU\}$. Then $A\sbst\BX$ is countable, and $\one\in\cl{A}$.
Thus, $\one\in\partlims(A)$.

As each $f\in A$ is $\{0,1\}$-valued, and limits of convergent sequences of $0$'s and $1$'s
must be either $0$ or $1$, each $f$ in $\partlims(A)$ is $\{0,1\}$-valued.
Let $C$ be the set of all partial $\{0,1\}$-valued functions $f$ on $X$, such that
$f\inv(1)\in\LI(\cU)$. Then $A\sbst C$, and $C$ is closed
under partial limits of sequences. Indeed, let $f_1,f_2,\dots\in C$, and $f=\lim_nf_n$.
As $\lim_nf_n(x)=f(x)$ and the functions $f_n$ are $\{0,1\}$-valued,
$f\inv(1)=\liminf_n f_n\inv(1)\in\LI(\cU)$.

Therefore, $\partlims(A)$ is contained in $C$, and in particular $\one\in C$,
that is, there is $B\in\LI(\cU)$ such that $X=\one\inv(1)\sbst B$. Thus, $X=B\in\LI(\cU)$.
\epf

\subsubsection*{Step 2: A selective property}

\bdfn
For a family $\cF$ of subsets of $X$, let
$$\cF_\downarrow=\{B\sbst X : (\exists A\in\cF)\ B\sbst A\},$$
the closure of $\cF$ under taking subsets.
\edfn

For a family $\cF$ of sets, $\Union\cF$ (without running index) denotes the union
of all members of $\cF$.
We say that a family of sets $\cV$ \emph{refines} another family $\cU$
if each $V\in\cV$ is contained in some $U\in\cU$. The following result
may be obtained by following arguments of Gerlits and Nagy \cite{GN}
and arguments of Nowik, Scheepers, and Weiss \cite{NSW}, proved for open
covers (under certain hypotheses on the space $X$).
We provide a different, direct proof, which makes no assumption on $X$.

\bprp\label{SP1}
Let $X$ be a topological space.
Assume that for each Borel $\omega$-cover $\cU$ of $X$, $X\in\LI(\cU)$.
Then for each sequence $\cU_1,\cU_2,\dots$ of Borel covers of $X$,
there are finite sets $\cF_1\sbst\cU_1,\cF_2\sbst\cU_2,\dots$, such that
for each $x\in X$, $x\in\Union\cF_n$ for all but finitely many $n$.
\eprp
\bpf
By moving to refinements, we may assume that for each $n$, the elements of $\cU_n$ are
pairwise disjoint, and $\cU_{n+1}$ refines $\cU_n$.\footnote{Given a Borel cover $\cU=\sseq{U_n}$,
the Borel cover $\sseq{U_n\sm(U_1\cup\dots\cup U_{n-1})}$ refines $\cU$, and its elements are pairwise
disjoint. Given two Borel covers $\cU,\cV$ whose elements are pairwise disjoint, the Borel
cover $\{U\cap V : U\in\cU, V\in\cV\}$ refines $\cU$ and $\cV$, and in particular its elements are pairwise
disjoint.}
This way, if there are infinitely many $n$ such that $\cU_n$ contains a finite subcover $\cF_n$ of $X$,
then this is true for all $n$ and the required assertion follows immediately.
Thus, we may assume that for each $n$, $\cU_n$ does not contain a finite subcover of $X$.

Let
$$\cB=\left\{\liminf_n\Union\cF_n : (\forall n)\ \cF_n\mbox{ is a finite subset of }\cU_n\right\}.$$
We must prove that $X\in\cB$.

$\LI(\cB_\downarrow)=\cB_\downarrow$:
For each $k$, assume that $B_k\sbst\liminf_n\Union\cF^k_n$, with each $\cF^k_n$ a finite subset
of $\cU_n$.
Take $\cF_n=\cF^1_n\cup\cF^2_n\cup\dots\cup\cF^n_n$ for each $n$. Then
$$\liminf_k B_k\sbst\liminf_k\liminf_n\Union\cF^k_n\sbst \liminf_n\Union\cF_n\in\cB,$$
and thus $\liminf_n B_n\in\cB_\downarrow$.

Thus, $\LI(\cB)\sbst\cB_\downarrow$, and therefore if $X\in\LI(\cB)$ then $X\in\cB$.
$\cB$ is an $\omega$-cover of $X$ and its elements are Borel,
but $\cB$ is in general \emph{not} countable, and thus we cannot apply the premise of the lemma.
To overcome this problem, we use a trick similar to one in \cite{GN}:
Define
$$\cA = \Union_{n\in\N}\left\{\Union\cF : \cF\sbst\cU_n, |\cF|=n\right\}.$$
$\cA$ is a Borel $\omega$-cover of $X$,
and therefore by the premise
of the lemma, $X\in\LI(\cA)\sbst\LI(\cA_\downarrow\cup\cB_\downarrow)$.
As $X\nin\cA$, it remains to show that
$\LI(\cA_\downarrow\cup\cB_\downarrow)=\cA_\downarrow\cup\cB_\downarrow$.

Let $B_1,B_2,\dots\in\cA_\downarrow\cup\cB_\downarrow$.
As $\liminf_n B_n\sbst\liminf_n B_{m_n}$ for each increasing
sequence $m_n$, and $\cA_\downarrow\cup\cB_\downarrow$ is closed downwards,
we may move to subsequences at our convenience.

If $B_n\in\cB_\downarrow$ for infinitely many $n$, then by moving to a subsequence
we may assume that $B_n\in\cB_\downarrow$ for all $n$, and therefore
$\liminf_n B_n\in\LI(\cB_\downarrow)=\cB_\downarrow\sbst\cA_\downarrow\cup\cB_\downarrow$.
In the remaining case, by moving to a subsequence, we may assume that $B_n\in\cA_\downarrow$ for all $n$.

Consider first the case where, after moving to an appropriate subsequence
of $B_1,B_2,\dots$, there is an increasing sequence $k_n$ such that $B_n\sbst\Union\cF_{k_n}$,
$\cF_{k_n}\sbst\cU_{k_n}$ with $|\cF_{k_n}|=k_n$, for all $n$.
As the covers $\cU_n$ are getting finer with $n$,
for each $i\notin\sseq{k_n}$ there is a finite $\cF_i\sbst\cU_i$
such that $\Union\cF_i$ contains $\Union\cF_{k_n}$ for the first
$n$ with $i<k_n$.
Then
$$\liminf_n B_n\sbst\liminf_n\Union\cF_n\in\cB,$$
as required.

Finally, there remains the case where, after moving to an appropriate subsequence
of $B_1,B_2,\dots$, there is $k$ such that for each $n$,
there is $\cF_n\sbst\cU_k$ with $|\cF_n|=k$, such that $B_n\sbst\Union\cF_n$.
Let $B=\lim\inf B_n$. We will show that $B\in\cA_\downarrow$.
We may assume that $B\neq\emptyset$. Take $x_1\in B$,
and $U_1\in\cU_k$ such that $x_1\in U_1$.
If $B\sbst U_1$, then $B\in\cA_\downarrow$. Otherwise, take $x_2\in B\sm U_1$,
and $U_2\in\cU_k$ such that $x_2\in U_2$.
Continue in the same manner until it is impossible to proceed,
but not more than $k$ steps, to have
$x_1,\dots,x_i\in B$, where $i\le k$, and distinct (and therefore disjoint) $U_1,\dots,U_i\in\cU_k$.
If $i<k$, then $B\sbst U_1\cup\dots\cup U_i$, a union of less than $k$ elements of $\cU_k$,
and thus $B\in\cA_\downarrow$.
Otherwise $i=k$, and for all but finitely many $n$, $x_1,\dots,x_k\in B_n\sbst\Union\cF_n$,
and as the elements of $\cU_k$ are pairwise disjoint, $\cF_n=\{U_1,\dots,U_k\}$ for all but finitely many
$n$. Consequently, $B\sbst\liminf_n\Union\cF_n= U_1\cup\dots\cup U_k\in\cA$, and therefore $B\in\cA_\downarrow$.
\epf

\subsubsection*{Step 3: A stronger selective property}
The selective property in the following theorem is stronger (\cite{CBC}, or Lemma \ref{tensor})
than the one introduced in the previous step.
In its original formulation \cite{MilProb}, Miller's Problem \ref{MP} asks whether the following
theorem is true.

\bthm\label{SP2}
Assume that for each Borel $\omega$-cover $\cU$ of $X$, $X\in\LI(\cU)$.
Then in fact, for each Borel $\omega$-cover $\cU$ of $X$, there
are $U_1,U_2,\dots\allowbreak\in\cU$ such that $X=\liminf_n U_n$.
\ethm
\bpf
Let
$$\cB=\{\liminf_n U_n : U_1,U_2,\dots\in\cU\}_\downarrow.$$
It suffices to show that $\LI(\cB)=\cB$.
Let $B_1,B_2,\dots\in\cB$, and $B=\liminf_n B_n$.
For each $n$, take $U^n_1,U^n_2,\dots\in\cU$ such that
$B_n\sbst\liminf_m U^n_m$.
Then for each $n$, the sets $V^n_m=\bigcap_{k\ge m} U^n_k$ are increasing to $B_n$,
and therefore the sets $V^n_m\cup (X\sm B_n)$ are increasing to $X$.

Applying Proposition \ref{SP1} to the covers $\cU_n=\{V^n_m\cup (X\sm B_n) : m\in\N\}$,
there are $m_n$ such that $X=\liminf_n V^n_{m_n}\cup (X\sm B_n)$
(since the covers are increasing, it suffices to pick one element from each cover).
As $\liminf_n B_n=B$, we have that
$$B\sbst (\liminf_n V^n_{m_n}\cup (X\sm B_n))\cap B\sbst \liminf_n V^n_{m_n}\sbst \liminf_n U^n_{m_n},$$
and therefore $B\in\cB$.
\epf

\subsubsection*{Step 4: Global to local}
The following lemma and its proof are, in the open/continuous case, due to Gerlits and Nagy \cite{GN}.
Their argument also applies to the Borel case.

\blem
Assume that for each Borel $\omega$-cover $\cU$ of $X$, there
are $U_1,U_2,\dots\in\cU$ such that $X=\liminf_n U_n$.
Then for each countable $A\sbst \BX$, each $f\in\cl{A}$ is a pointwise limit of a sequence
of elements of $A$.
\elem
\bpf
We may assume, by adding the function $\one-f$ to all considered functions,
that $f=\one$, the constant $1$ function.
For each $n$, let $\cU_n=\{g\inv[(1-1/n,1+1/n)] : g\in A\}$.
As $\one\in\cl{A}$, $\cU_n$ is a (Borel) $\omega$-cover of $X$.
By Theorem \ref{SP2}, there are $g_n\in A$ such that $X=\liminf_n g_n\inv[(1-1/n,1+1/n)]$.
Then $\one=\lim_n g_n$.
\epf

This completes the proof of Theorem \ref{main}.

\section{Baire class 1 functions (Herna\'ndez's Problem)}\label{HernSol}

The following Theorem, which strengthens Theorem \ref{main} (in the realm of perfectly normal spaces),
answers in the positive a question of Salvador Hern\'andez.

A topological space $X$ is \emph{perfectly normal} if it is normal (any two disjoint
closed sets have disjoint neighborhoods), and each open subset of $X$ is $F_\sigma$,
that is, a union of countably many closed subsets of $X$.
For example, metric spaces are perfectly normal.

A function $f:X\to\R$ is of \emph{Baire class 1} if $f$ is the pointwise limit of
a sequence of continuous real-valued functions on $X$.
Let $\BaireOneX\sbst\R^X$ denote the subspace of all Baire class 1 functions $f:X\to\R$.

\bthm\label{B1}
Let $X$ be a perfectly normal topological space.
Assume that, for each countable $A\sbst\BaireOneX$, each $f\in\cl{A}$ (closure in $\BaireOneX$)
belongs to the closure of $A$ under partial limits of sequences.
Then for each countable $A\sbst \BaireOneX$, each $f\in\cl{A}$ (closure in $\BaireOneX$)
is a limit of a sequence of elements of $A$.

Moreover, for each countable $A\sbst \BX$, each $f\in\cl{A}$ (closure in $\BX$) is a limit of a sequence
of elements of $A$.
\ethm
\bpf
As the closure in a subspace $Y$ of $\R^X$ is equal to the intersection of the closure in $\R^X$ and $Y$,
and $\BaireOneX\sbst \BX$, it suffices to prove the second assertion. We follow the proof steps of Theorem \ref{main},
and modify them when needed.

A set $A\sbst X$ is $\Delta^0_2$ if both $A$ and $X\sm A$ are $F_\sigma$.
The family $\Delta^0_2(X)$ of all $\Delta^0_2$ subsets of $X$ forms an algebra
of sets, that is, it is closed under finite unions and complements (and therefore
also under finite intersections and set differences). This fact is
applied repeatedly when following the steps in the proof of Theorem \ref{main}.

A function $f:X\to\R$ is \emph{$\Delta^0_2$-measurable} if for each open $U\sbst\R$,
$f\inv[U]$ is $\Delta^0_2$.
For each $\Delta^0_2$ set $U\sbst X$, $\chi_U$ is $\Delta^0_2$-measurable.

The following lemma is proved for the metrizable case in \cite[Lemma 24.12]{Kechris}.
The proof there uses only Urysohn's Lemma, which applies for all normal spaces.

\blem[folklore]\label{Bc1}
Let $X$ be a normal space, and $U$ be a $\Delta^0_2$ subset of $X$.
Then $\chi_U$ is of Baire class 1.
\elem
\bpf
Let $F_n\sbst X$ be closed, and $G_n\sbst X$ be open, such that $F_n\sbst F_{n+1}\sbst U\sbst G_{n+1}\sbst G_n$ for all $n$,
and $U=\Union_nF_n=\bigcap_n G_n$.
By Urysohn's Lemma, there is for each $n$ a continuous function $f_n:X\to\R$ such that
$f_n(x)=1$ for all $x\in F_n$ and $f_n(x)=0$ for all $x\nin G_n$.
Then $\lim_n f_n(x)=\chi_U(x)$ for all $x\in X$.
\epf

Thus, arguing as in Step 1 of Theorem \ref{main}, we have that
for each countable $\Delta^0_2$ $\omega$-cover $\cU$ of $X$, $X\in\LI(\cU)$.

The arguments of Step 2 show the following.

\bprp\label{D02}
Assume that for each countable $\Delta^0_2$ $\omega$-cover $\cU$ of $X$, $X\in\LI(\cU)$.
Then for each sequence $\cU_1,\cU_2,\dots$ of countable $\Delta^0_2$ covers of $X$,
there are finite sets $\cF_1\sbst\cU_1,\cF_2\sbst\cU_2,\dots$, such that
for each $x\in X$, $x\in\Union\cF_n$ for all but finitely many $n$.\qed
\eprp

In particular, as $X$ is perfectly normal, $X$ has the property in the conclusion of Proposition \ref{D02}
for \emph{closed} sets. We use the following strong result of Bukovsk\'y, Rec\l{}aw, and Repick\'y.

\blem[Bukovsk\'y--Rec\l{}aw--Repick\'y \cite{BRR91}]\label{surprise1}
Let $X$ be a perfectly normal space.
Assume that for each sequence $\cU_1,\cU_2,\dots$ of countable closed covers of $X$,
there are finite sets $\cF_1\sbst\cU_1,\cF_2\sbst\cU_2,\dots$, such that
for each $x\in X$, $x\in\Union\cF_n$ for all but finitely many $n$.
Then the same holds for each sequence $\cU_1,\cU_2,\dots$ of Borel covers of $X$.
\elem

The property established in Lemma \ref{surprise1} implies that every Borel subset of
$X$ is $F_\sigma$ (e.g., \cite{CBC}), and thus every Borel set is $\Delta^0_2$.
By the property established before Proposition \ref{D02}, we have that,
for each countable \emph{Borel} $\omega$-cover $\cU$ of $X$, $X\in\LI(\cU)$.

Thus, theorem \ref{SP2} and Step 4 apply, and the proof is completed.
\epf

\brem
The proof of Theorem \ref{B1} shows that it suffices to assume
that for each countable set $A$ of $\Delta^0_2$-measurable real-valued functions
on $X$, the closure of $A$ in the space of all $\Delta^0_2$-measurable real-valued functions
on $X$ is contained in $\partlims(A)$.
\erem

\section{Continuous functions (Gerlits--Nagy's Problem)}\label{OnGNP}

Thus far, we have refrained from using the notation of the field of selective properties,
despite their playing important role in the proofs.
However, as we are about to make a more extensive use of the theory,
we give here the necessary introduction.
Readers who wish to learn more on the topic and its history are referred to any of its surveys
\cite{LecceSurvey, KocSurv, ict}.

Let $X$ be a topological space. 
Let $\rmO(X)$ be the family of all open covers of $X$.
Define the following subfamilies of $\rmO(X)$: $\cU\in\Om(X)$ if $\cU$ is an $\omega$-cover of $X$.
$\cU\in\Ga(X)$ if $\cU$ is infinite, and each element of $X$ is contained in all but finitely many members of $\cU$.

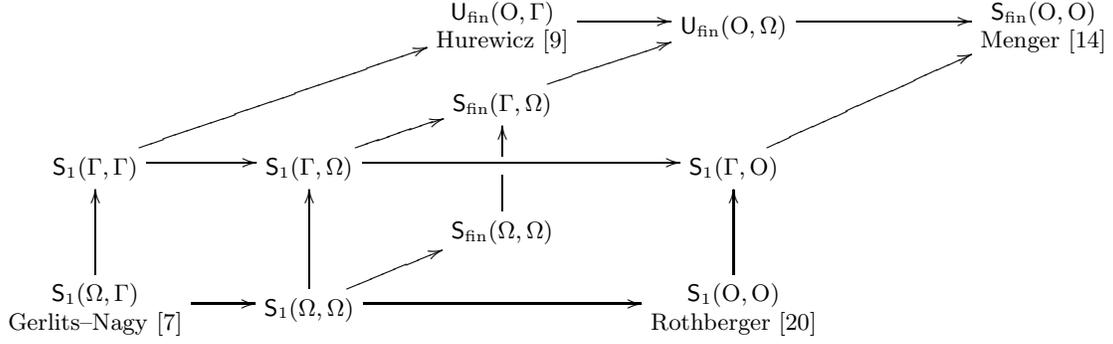
\begin{figure}[!htp]
\begin{changemargin}{-4cm}{-3cm}
\begin{center}
{\scriptsize
$\xymatrix@R=8pt{
&
&
& \sr{\ufin(\rmO,\Ga)}{Hurewicz \cite{Hure25}}\ar[r]
& \sr{\ufin(\rmO,\Om)}{}\ar[rr]
& & \sr{\sfin(\rmO,\rmO)}{Menger \cite{Menger24}}
\\
&
&
& \sr{\sfin(\Ga,\Om)}{}\ar[ur]
\\
& \sr{\sone(\Ga,\Ga)}{}\ar[r]\ar[uurr]
& \sr{\sone(\Ga,\Om)}{}\ar[rr]\ar[ur]
& & \sr{\sone(\Ga,\rmO)}{}\ar[uurr]
\\
&
&
& \sr{\sfin(\Om,\Om)}{}\ar'[u][uu]
\\
& \sr{\sone(\Om,\Ga)}{Gerlits--Nagy \cite{GN}}\ar[r]\ar[uu]
& \sr{\sone(\Om,\Om)}{}\ar[uu]\ar[rr]\ar[ur]
& & \sr{\sone(\rmO,\rmO)}{Rothberger \cite{Roth41}}\ar[uu]
}$
}
\caption{The Scheepers Diagram}\label{SchDiag}
\end{center}
\end{changemargin}
\end{figure}

Some of the following statements may hold for families $\scrA$ and $\scrB$ of covers of $X$.
\begin{description}
\item[$\binom{\scrA}{\scrB}$] Each element of $\scrA$ contains an element of $\scrB$.
\medskip
\item[$\sone(\scrA,\scrB)$] For all $\cU_1,\cU_2,\dots\in\scrA$, there are
$U_1\in\cU_1,U_2\in\cU_2,\dots$ such that $\sseq{U_n}\in\scrB$.
\medskip
\item[$\sfin(\scrA,\scrB)$] For all $\cU_1,\cU_2,\dots\in\scrA$, there are
finite $\cF_1\sbst\cU_1,\cF_2\sbst\cU_2,\dots$ such that $\Union_n\cF_n\in\scrB$.
\medskip
\item[$\ufin(\scrA,\scrB)$] For all $\cU_1,\cU_2,\dots\in\scrA$, none containing
a finite subcover, there are finite $\cF_1\sbst\cU_1,\cF_2\sbst\cU_2,\dots$ such that $\sseq{\Union\cF_n}\in\scrB$.
\end{description}
We say, e.g., that $X$ satisfies $\sone(\rmO,\rmO)$ if the statement $\sone(\rmO(X),\rmO(X))$ holds.
This way, $\sone(\rmO,\rmO)$ is a property of topological spaces, and similarly for all other statements
and families of covers.
Under some mild hypotheses on the considered topological spaces,
each nontrivial property among these properties, where $\scrA,\scrB$ range over $\rmO,\Om,\Ga$,
is equivalent to one in Figure \ref{SchDiag}, named after Scheepers in recognition
of his seminal contribution to the field. In this diagram, an arrow denotes implication.

Other types of covers, most notably \emph{Borel} covers, were also considered in this
context. We say, for example, that $X$ satisfies $\sone(\Om,\Om)$ \emph{for Borel covers}
if $\sone(\Om(X),\Om(X))$ holds, when redefining $\Om(X)$ to consist of all countable 
Borel $\omega$-covers of $X$.

For clarity of notation, we identify a property with the family of topological spaces (of a certain type, which should be clear
from the context) satisfying it.

The property deduced in Theorem \ref{SP1} is $\ufin(\rmO,\Ga)$ for Borel covers.
For Borel covers, $\ufin(\rmO,\Ga)=\sone(\Ga,\Ga)$ \cite{CBC}, and using this
the proof of Theorem \ref{SP2} can be slightly simplified.

Gerlits and Nagy \cite{GN} proved the following lemma for Hausdorff spaces.
We will see that it holds for arbitrary topological spaces.

\blem\label{tensor}
$\binom{\Om}{\Ga}=\sone(\Om,\Ga)$ (for general topological spaces).
\elem
\bpf
Assume that $X$ satisfies $\binom{\Om}{\Ga}$, and let $\cU_1,\cU_2,\dots$ be open $\omega$-covers of $X$.
We may assume that for each $n$, $\cU_{n+1}$ refines $\cU_n$.

For each $n$, enumerate $\cU_n=\{U^n_m : m\in\N\}$.
Let $V_m=U^1_m$ for all $m$.
Define
$$\cW=\Union_{n\in\N}\{V_n\cap U^n_m : m\in\N\}.$$
$\cW$ is an open $\omega$-cover of $X$.
Thus, there are $W_1,W_2,\dots\in\cW$ such that $X=\liminf_k W_k$.
Fix $n$. As $V_n\neq X$, it is not possible that $W_k\in \{V_n\cap U^n_m : m\in\N\}$
for infinitely many $k$. Since the sets $U^n_m$ are increasing with $m$,
we may assume that there is at most one $W_k$ in each set $\cW_n=\{V_n\cap U^n_m : m\in\N\}$.
For each $n$, let $r_n\ge n$ be the first such that there is some $W_k$ in $\cW_{r_n}$.
Since the covers $\cU_n$ get finer with $n$, we can pick for each $n$ an
element $U^n_{m_n}\in\cU_n$ containing the $W_k$ which is in $\cW_{r_n}$.
Then $X=\liminf_kW_k\sbst\liminf_n U^n_{m_n}$, and therefore $\liminf_n U^n_{m_n}=X$.
\epf

Using Lemma \ref{tensor}, Gerlits and Nagy proved the following fundamental local-to-global correspondence result.

\bthm[Gerlits--Nagy \cite{GN}]\label{1}
For Tychonoff spaces $X$, the following properties are equivalent:
\be
\itm For each $A\sbst \CX$, each $f\in\cl{A}$ is a limit of a sequence of elements of $A$ (i.e.,
$\CX$ is Fr\'echet--Urysohn).
\itm $X$ satisfies $\binom{\Om}{\Ga}$.
\ee
\ethm

We establish a similar result for the other major property studied in the present paper.
To this end, we need the following definition and a lemma.

\bdfn
$\rmL(X)$ is the family of open covers of $X$ such that $X\in\LI(\cU)$.
\edfn

Theorem \ref{SP2} tells that $\binom{\Om}{\rmL}=\binom{\Om}{\Ga}$ for Borel covers.
In particular, using that $\binom{\Om}{\Ga}=\sone(\Om,\Ga)$ for Borel covers,
we have that $\binom{\Om}{\rmL}=\sone(\Om,\rmL)$ for Borel covers.
The last assertion also holds in the open case, but a different proof is required.

\blem\label{sel1}
$\binom{\Om}{\rmL}=\sone(\Om,\rmL)=\sfin(\Om,\rmL)$.
\elem
\bpf
As $\sone(\Om,\rmL)$ implies $\sfin(\Om,\rmL)$, which in turn implies $\binom{\Om}{\rmL}$,
it remains to prove that $\binom{\Om}{\rmL}$ implies $\sone(\Om,\rmL)$.
To this end, it suffices to prove that $\binom{\Om}{\rmL}$ implies $\sone(\Om,\Om)$.
$\sone(\Om,\Om)$ is equivalent to having all finite powers of $X$ satisfy $\sone(\rmO,\rmO)$ \cite{Sakai88}.
Gerlits and Nagy \cite{GN} proved that $\binom{\Om}{\rmL}$ implies $\sone(\rmO,\rmO)$.
Thus, it remains to prove that $\binom{\Om}{\rmL}$ is preserved by finite powers.

Assume that $X$ satisfies $\binom{\Om}{\rmL}$, and let $k\in\N$.
Let $\cU$ be an open $\omega$-cover of $X^k$. Then there is an open
$\omega$-cover $\cV$ of $X$ such that $\cV'=\{V^k : V\in\cV\}$ refines $\cU$ \cite{coc2}.
Then $X\in\LI(\cV)$. For arbitrary sets $B_1,B_2,\dots$, $\liminf_n (B_n)^k=(\liminf_n B_n)^k$.
Thus, $X^k\in\{B^k : B\in\LI(\cV)\}=\LI(\cV')$, and therefore $X^k\in\LI(\cU)$.
\epf

\bthm\label{2}
For Tychonoff spaces $X$, the following properties are equivalent:
\be
\itm For each $A\sbst \CX$, $\cl{A}\sbst\partlims(A)$.
\itm $X$ satisfies $\binom{\Om}{\rmL}$ (that is, for each open $\omega$-cover $\cU$ of $X$, $X\in\LI(\cU)$).
\ee
\ethm
\bpf
$(1\Impl 2)$ For partial functions $f$ and $g$, $g\circ f$ is the partial
function with domain $\{x\in\dom(f) : f(x)\in\dom(g)\}$, defined as
usual by $g\circ f(x)=g(f(x))$.

For a surjection $\vphi:X\to Y$ and partial functions $f_n:Y\to\R$,
the domain of $\lim_n(f_n\circ \vphi)$ is $\vphi\inv[\dom(\lim_nf_n)]$,
and $\lim_n(f_n\circ \vphi)=(\lim_nf_n)\circ\vphi$. Thus, we have the following.

\blem\label{c}
Assume that for each $A\sbst \CX$, each $\cl{A}\sbst\partlims(A)$.
Then every continuous image of $X$ has the same property.\qed
\elem

A topological space is \emph{zero-dimensional} if
its clopen (simultaneously closed and open) sets form a base for its topology.
An argument similar to one in \cite{GN} gives the following.

\blem
Let $X$ be a Tychonoff space.
Assume that for each $A\sbst \CX$, each $f\in\cl{A}$ belongs to the
closure of $A$ under partial limits of sequences.
Then $X$ is zero-dimensional.
\elem
\bpf
It suffices to prove that $[0,1]$ is not a continuous image of $X$.
Indeed, for each open $U\sbst X$ and each $a\in U$, let $\Psi:X\to [0,1]$ be
continuous, such that $\Psi(a)=0$ and $\Psi(x)=1$ for all $x\in X\sm U$.
Take $r\in [0,1]$ which is not in the image of $\Psi$.
Then $\Psi\inv[[0,r)]$ is a clopen neighborhood of $x$ contained in $U$.

Assume that $[0,1]$ is a continuous image of $X$.
Let $A\sbst\op{C}([0,1])$ be the set of all continuous
$f:[0,1]\to[0,1]$ such that the Lebesgue measure of $f\inv[(1/2,1]]$ is at most $1/2$.
Then $\one$ is in the closure of $A$.
Let $C$ be the set of all \emph{partial} $f:[0,1]\to[0,1]$ such that $f\inv[(1/2,1]]$ is Lebesgue measurable,
and its measure is at most $1/2$. $C$ is closed under partial limits of sequences and contains $A$,
but $\one\nin C$; a contradiction.
\epf

Let $\cU$ be an open $\omega$-cover of $X$. As $X$ is zero-dimensional, $\cU$ can be refined
to a clopen $\omega$-cover of $X$ by replacing each $U\in\cU$ with all finite unions of clopen subsets of $U$.
Now, for each clopen $U$ the function $\chi_U$ is continuous,
and $\one$ is in the closure of $\{\chi_U : U\in\cU\}$. By (1), $\one$ is in the closure
of $\{\chi_U : U\in\cU\}$ under partial limits of sequences.
Continue as in the proof of Lemma \ref{l2g}.

$(2\Impl 1)$ In (1), by adding $\one-f$ to all of the involved partial functions,
it suffices to consider the case $f=\one$.
Let $A\sbst \CX$, and assume that $\one\in\cl{A}$.
For each $n$, let $\cU_n=\{f\inv[(1-1/n,1+1/n)] : f\in A\}$.
$\cU_n$ is an open $\omega$-cover of $X$. By Lemma \ref{sel1},
there are $f_1,f_2,\dots\in A$ such that $X\in\LI(\sseq{f_n\inv[(1-1/n,1+1/n)]})$.

We claim that
$$\cA=(\sseq{f_n\inv[(1-1/n,1+1/n)]}\cup\{f\inv(1) : f\in\partlims(A)\})_\downarrow$$
is closed under the operator $\liminf$. Indeed,
assume that we are given a sequence of elements of $\cA$. By thinning it out,
and replacing each element by an appropriate element containing it,
we may assume that this sequence is all in $\sseq{f_n\inv[(1-1/n,1+1/n)]}$ or all in $\{f\inv(1) : f\in\partlims(A)\}$.
In the first case, by thinning out further we may assume that the sequence is either constant (in which
case we are done), or consists of distinct elements $f_{m_n}\inv[(1-1/{m_n},1+1/{m_n})]$ with $m_n$ increasing.
In this case, let $f=\lim_nf_{m_n}$. For each $x\in\liminf_n f_{m_n}\inv[(1-1/{m_n},1+1/{m_n})]$,
$f(x)=\lim_nf_{m_n}(x)=1$, and thus $\liminf_n f_{m_n}\inv[(1-1/{m_n},1+1/{m_n})]$ is in $\cA$.
The second case is similar (and slightly easier).

Thus, $X\in\cA$, which means that there is $f\in\partlims(A)$ such that $X=f\inv(1)$, that is,
$\one=f\in\partlims(A)$.
\epf

Clearly, $\binom{\Om}{\Ga}$ implies $\binom{\Om}{\rmL}$.
The original Gerlits--Nagy Problem, posed in \cite{GN}, asks whether
these properties are in fact equivalent (for Tychonoff $X$, or even for $X\sbst\R$).
Theorems \ref{1} and \ref{2} justify the reformulation given in Problem \ref{GNP}.

Originally, Gerlits and Nagy \cite{GN} studied five properties, numbered $\alpha,\beta,\gamma,\delta,\epsilon$,
where each property implies the subsequent one. $\binom{\Om}{\Ga}$ and $\binom{\Om}{\rmL}$ were
numbered $\gamma$ and $\delta$, respectively, and are often named accordingly in the literature.
Their problem was originally stated as whether property $\delta$ implies (and is therefore equivalent to)
property $\gamma$.

A topological space $X$ is said to satisfy a property $P$ \emph{hereditarily} if each $Y\sbst X$ satisfies $P$.
Pushing our methods further, we can solve the Gerlits--Nagy Problem in the affirmative
for spaces $X$ satisfying $\sone(\Ga,\Ga)$ hereditarily. We will use the following
result of Francis Jordan \cite{Jordan} (see also \cite{LinSAdd}), proved using a
new fusion argument of his.

\blem[Jordan]\label{jlem}
Let $B=\Union_nB_n\sbst X$ be an increasing union, where each $B_n$ satisfies
$\sone(\Gamma,\Gamma)$. For all open sets $U^n_m\sbst X$, $n,m\in\N$,
with $B_n\sbst\liminf_m U^n_m$ for each $n$, there are $m_1,m_2,\dots\in\N$
such that $B\sbst\liminf_n U^n_{m_n}$.
\elem

\bthm\label{wGNP}
For topological spaces $X$ satisfying $\sone(\Ga,\Ga)$ hereditarily, the following are equivalent:
\be
\itm $X$ satisfies $\binom{\Om}{\rmL}$.
\itm $X$ satisfies $\binom{\Om}{\Ga}$.
\ee
\ethm
\bpf[Proof of $(1)\Impl(2)$]
\blem\label{LG}
Assume that $X$ satisfies $\sone(\Ga,\Ga)$ hereditarily. Then $X$ satisfies $\binom{\rmL}{\Ga}$.
\elem
\bpf
Let $\cU$ be an open cover of $X$ with $X\in\LI(\cU)$. Define
$$\cB=\{\liminf_n U_n : U_1,U_2,\dots\in\cU\}.$$
We will prove that $X\in\cB$. To this end, it suffices to show that $\LI(\cB_\downarrow)=\cB_\downarrow$.

Let $B_1,B_2,\dots\in\cB_\downarrow$, and $B=\liminf_n B_n$.
Replacing each $B_n$ with $\bigcap_{m\ge n}B_m$, we may assume
that $B_1\sbst B_2\sbst \dots$, and $\Union_n B_n=B$.
For each $n$, take $U^n_1,U^n_2,\dots\in\cU$ such that
$B_n\sbst\liminf_m U^n_m$. By the premise of the proposition,
each $B_n$ satisfies $\sone(\Ga,\Ga)$.
By Jordan's Lemma \ref{jlem}, there are $m_1,m_2,\dots\in\N$
such that $B\sbst_n\liminf_n U^n_{m_n}\in\cB_\downarrow$, and therefore $B\in\cB_\downarrow$.
\epf
It remains to note that the conjunction of $\binom{\rmL}{\Ga}$ and $\binom{\Om}{\rmL}$ implies $\binom{\Om}{\Ga}$.
\epf

\brem
For each topological space $X$, $\Ga(X)\sbst\rmL(X)\sbst\Om(X)$. To see the second inclusion,
assume that there is a finite $F\sbst X$ not covered by any $U\in\cU$. Then $F$ is not covered by
any element of $\LI(\cU)$, and in particular, $X\notin \LI(\cU)$.
Thus, the implication at the end of the proof of Theorem \ref{wGNP} is in fact
an equivalence, that is, $\binom{\Om}{\rmL}\cap\binom{\rmL}{\Ga}=\binom{\Om}{\Ga}$.
\erem

\bcor\label{2h}
For Tychonoff spaces $X$ satisfying $\sone(\Ga,\Ga)$, the following are equivalent:
\be
\itm $X$ satisfies $\binom{\Om}{\rmL}$ hereditarily.
\itm $X$ satisfies $\binom{\Om}{\Ga}$ hereditarily.
\ee
\ecor
\bpf[Proof of $(1)\Impl(2)$]
By Theorem \ref{wGNP}, it suffices to prove that $X$ satisfies $\sone(\Ga,\Ga)$ \emph{hereditarily}.

Nowik, Scheepers and Weiss  proved that $\binom{\Om}{\rmL}$ implies $\ufin(\rmO,\Ga)$ \cite{NSW}.\footnote{For
a direct proof, see the proof of Proposition \ref{SP1}.}
Thus, if $X$ satisfies $\binom{\Om}{\rmL}$ hereditarily, then $X$ satisfies $\ufin(\rmO,\Ga)$ hereditarily.
Fremlin and Miller \cite{FM88} proved that in the latter case, $X$ is a $\sigma$-space, that is, each Borel subset of
$X$ is $F_\sigma$. This, together with $X$'s satisfying $\sone(\Ga,\Ga)$, implies that $X$ satisfies $\sone(\Ga,\Ga)$ hereditarily
\cite{Hales05, LinSAdd}.
\epf

\brem
The argument in the proof of Corollary \ref{2h} shows that, for Tychonoff \emph{$\sigma$-spaces} $X$,
$\binom{\Om}{\Ga}=\binom{\Om}{\rmL}\cap\sone(\Ga,\Ga)$. In this case, this joint property
coincides with its hereditary version.
\erem

Assuming that the answer to the Gerlits--Nagy Problem is \emph{negative},
the results of this section explain, to some extent, why no counter example was discovered thus far.
A natural strategy would be to begin with a set
$X\sbst\R$ satisfying $\binom{\Om}{\Ga}$, and then look for a subset of $X$,
in a way which ``destroys'' $\binom{\Om}{\Ga}$, but not too much, so that
$\binom{\Om}{\rmL}$ still holds.
There are several constructions of subsets of $\R$ satisfying $\binom{\Om}{\Ga}$.
The first one is due to Galvin and Miller \cite{GM}. Here, $X$ has a countable subset $Q$ such that $X\sm Q$ does not satisfy $\binom{\Om}{\Ga}$.
Unfortunately, $X\sm Q$ does not even satisfy $\ufin(\rmO,\Ga)$, and in particular not $\binom{\Om}{\rmL}$.\footnote{On the other hand,
we proved in \cite{LinSAdd} that any ``natural'' change of Galvin and Miller's construction \emph{without} moving
to a subset at the end would keep $X$ in $\binom{\Om}{\Ga}$.}
Another, substantially different, construction is due to Todor\v{c}evic \cite{GM}, but this $X$ satisfies $\binom{\Om}{\Ga}$
hereditarily. Finally, using a variation of Todor\v{c}evic's method,
Miller \cite{MilNonGamma} constructed $X\sbst\R$ satisfying $\binom{\Om}{\Ga}$ for Borel covers,
and a subset $Y$ of $X$ not satisfying $\binom{\Om}{\Ga}$.
$\binom{\Om}{\Ga}$ for Borel covers, implies $\sone(\Ga,\Ga)$ for Borel covers, which is hereditary.
Thus, in this case $X$ satisfies $\sone(\Ga,\Ga)$ \emph{hereditarily},
and by Theorem \ref{wGNP} no subset of $X$ would separate $\binom{\Om}{\rmL}$ from $\binom{\Om}{\Ga}$.

We conclude this section with a local reformulation of Theorem \ref{wGNP}.
A topological space $Z$ has the Arhangel'ski\u{\i} property $\alpha_2$ if, for each $z\in Z$,
whenever $\lim_m z^n_m=z$ for all $n$, there are $m_1,m_2,\dots$ such that $\lim_n z^n_{m_n}=z$.
When $Z=\CX$, we can take $z=\one$ in the definition. Hale\v{s} proved that, for perfectly normal spaces
$X$, the following properties are equivalent:
\be
\itm For each $Y\sbst X$, $\CY$ is an $\alpha_2$ space.
\itm $X$ satisfies $\sone(\Ga,\Ga)$ hereditarily.
\ee
Collecting together the results of this section, we have the following.

\bthm
Let $X$ be a perfectly normal space, such that for each $Y\sbst X$, $\CY$ is an $\alpha_2$ space.
Then the following properties are equivalent:
\be
\itm For each $A\sbst \CX$, $\cl{A}\sbst\partlims(A)$.
\itm For each $A\sbst \CX$, each $f\in\cl{A}$ is a limit of a sequence of elements of $A$ (i.e., $\CX$ is
Fr\'echet--Urysohn).\qed
\ee
\ethm

\section{Some results about the missing piece}

The property $\binom{\rmL}{\Ga}$ was central, implicitly or explicitly, in our proofs,
for the basic reason that $$\binom{\Om}{\Ga}=\binom{\Om}{\rmL}\cap\binom{\rmL}{\Ga}.$$
To prove that $\binom{\Om}{\Ga}=\binom{\Om}{\rmL}$ (the Gerlits--Nagy Problem),
it is necessary and sufficient to prove that $\binom{\Om}{\rmL}$ implies $\binom{\rmL}{\Ga}$.
We therefore describe some fundamental properties of  $\binom{\rmL}{\Ga}$, and the ensuing
open problems concerning it.

\bprp\label{qs}
$\binom{\rmL}{\Ga}=\sone(\rmL,\Ga)=\sfin(\rmL,\Ga)$. In particular, $\binom{\rmL}{\Ga}$ implies $\sone(\Ga,\Ga)$.
\eprp
\bpf
It suffices to prove the last assertion.
Assume that for each $n$, $\cU_n=\{U^n_m : m\in\N\}\in\Ga(X)$.
We may assume that the covers $\cU_n$ get finer with $n$.\footnote{If $\sseq{U_n},\sseq{V_n}\in\Ga(X)$,
then $\sseq{U_n\cap V_n}\in\Ga(X)$ and is finer than both.}

Let $V_m=U^1_m$ for all $m$. Define
$$\cW=\Union_{n\in\N}\{V_n\cap U^n_m : m\in\N\}.$$
Then
$$\liminf_n\liminf_m (V_n\cap U^n_m) = \liminf_n V_n=X,$$
and therefore $X\in\LI(\cW)$. By $\binom{\rmL}{\Ga}$,
there are $W_1,W_2,\dots\in\cU$ such that $X=\liminf_k W_k$.
Fix $n$. As $V_n\neq X$, it is not possible that $W_k\in \{V_n\cap U^n_m : m\in\N\}$
for infinitely many $k$. Thus, by thinning out the sequence $W_k$ if needed, we may
assume that there is at most one $W_k$ in each set $\{V_n\cap U^n_m : m\in\N\}$.
Since the covers $\cU_n$ get finer with $n$, we can pick for each $n$ an
element $U^n_{m_n}\in\cU_n$, such that $X=\liminf_n U^n_{m_n}$.
\epf

\bprp\label{nss}
The property of satisfying $\sone(\Ga,\Ga)$ hereditarily is strictly stronger than $\binom{\rmL}{\Ga}$.
\eprp
\bpf
Lemma \ref{LG} tells that hereditarily-$\sone(\Ga,\Ga)$ implies $\binom{\rmL}{\Ga}$.
Assuming for example \CH{}, there is $X\sbst\R$ and a subset $Y$ of $X$ such that $X$ satisfies $\binom{\Om}{\Ga}$
(and thus also $\binom{\rmL}{\Ga}$), and $Y$ does not even satisfy $\sfin(\rmO,\rmO)$, and in particular not $\sone(\Ga,\Ga)$
\cite{GM}. Apply Proposition \ref{qs}.
\epf

If $\sone(\Ga,\Ga)$ implies $\binom{\rmL}{\Ga}$, then the word ``hereditarily'' can be removed from Theorem \ref{wGNP}.
However, we suspect that this is not the case.

\bcnj
$\binom{\rmL}{\Ga}$ is strictly stronger than $\sone(\Ga,\Ga)$.
\ecnj

To prove this conjecture, it suffices to construct (say using \CH) sets $X,Y\sbst\R$ satisfying
$\binom{\rmL}{\Ga}$, such that $X\cup Y$ does not satisfy $\binom{\rmL}{\Ga}$, because $\sone(\Ga,\Ga)$
is $\sigma$-additive.

\bprb
Is $\binom{\rmL}{\Ga}$ preserved by finite unions?
\eprb

If it is, then $\binom{\rmL}{\Ga}$ is in fact $\sigma$-additive, because of the following.

\bprp
$\binom{\rmL}{\Ga}$ is linearly $\sigma$-additive, that is, is preserved by countable increasing unions.
\eprp
\bpf
Assume that $X_1\sbst X_2\sbst\dots$ all satisfy $\binom{\rmL}{\Ga}$, $\Union_nX_n=X$, and $X\in\LI(\cU)$.
Then for each $n$, $X_n\in\LI(\{U\cap X_n  : U\in\cU\})$, and thus there
are $U^n_m\in\cU$, $m\in\N$, such that $X_n\sbst\liminf_m U^n_m$.
By Jordan's Lemma \ref{jlem}, there are $m_1,m_2,\dots\in\N$ such that
$X=\liminf_n U^n_{m_n}$.
\epf

The proofs of the above results are also valid in the case of Borel covers,
and since $\sone(\Ga,\Ga)$ for Borel covers is hereditary, we have the following.

\bcor
For Borel covers, $\binom{\rmL}{\Ga}=\sone(\Ga,\Ga)$.\qed
\ecor

Thus, none of the above-mentioned problems remains open in the Borel case.

We conclude with a local characterization of $\binom{\rmL}{\Ga}$.

\bthm
For perfectly normal spaces $X$, the following are equivalent.
\be
\itm For each $A\sbst\CX$ each $f\in\CX\cap\partlims(A)$ is a limit of a sequence of elements of $A$.
\itm $X$ satisfies $\binom{\rmL}{\Ga}$.
\ee
\ethm
\bpf
$(1\Impl 2)$

\blem\label{cl}
Let $X$ be a perfectly normal space.
Assume that for each $A\sbst\CX$, each $f\in\CX\cap\partlims(A)$ is a limit of a sequence of elements of $A$.
Then each element of $\rmL(X)$ has a clopen refinement in $\rmL(X)$.
\elem
\bpf
Indeed, this follows from a formally weaker property: Let $P$ be the property that,
for each $A\sbst\CX$, each $f$ in the closure of $A$ \emph{in $\CX$} under limits of sequences,
is a limit of a sequence of elements of $A$.

Fremlin \cite{FrwQN} proved that
$P$ is equivalent to the property named wQN in \cite{BRR91}, where it is shown that for perfectly
normal spaces, wQN implies that each open set is a countable union of \emph{clopen} sets \cite[Corollary 4.6]{BRR91}.

Now, let $\cU\in\rmL(X)$. For each $U\in\cU$, present $U$ as an increasing union $U=\Union_nC_n(U)$ of clopen sets.
Then $U=\liminf_n C_n(U)$. Let $\cV=\{C_n(U) : U\in\cU, n\in\N\}$. Then $\cV$ is a clopen refinement of $\cU$, and
$X\in\LI(\cU)\sbst\LI(\cV)$, that is, $\cV\in\rmL(X)$.
\epf

Let $\cU\in\rmL(X)$. By Lemma \ref{cl}, we may assume that the elements of $\cU$ are clopen.
Let $A=\{\chi_U : U\in\cU\}$. $A\sbst\CX$.
Let $\cV=\{f\inv(1) : f\in\partlims(A)\}$. $\cU\sbst\cV$, and $\cV$ is closed under the operator $\liminf$.
Indeed, Let $f_1,f_2,\dots\in\partlims(A)$, and $B=\liminf_n f_n\inv(1)$.
As $f=\lim_nf_n\in\partlims(A)$, $B=f\inv(1)\in\cV$.

Thus, $X\in\cV$, and therefore $\one\in\partlims(A)$.
By (1), there are $U_n\in\cU$ such that $\lim_n\chi_{U_n}=\one$, that is, $\liminf_n U_n=X$.

$(2\Impl 1)$
Assume that $\one\in\partlims(A)$.
For each $n$, let $\cU_n=\{f\inv[(1-1/n,1+1/n)] : f\in A\}$. $\cU_n\in\rmL(X)$.
Indeed, let $C$ be the family of all partial $f:X\to\R$, such that $f\inv[(1-1/n,1+1/n)]\in\LI(\cU_n)$.
Then $A\sbst C$, and $C$ is closed under partial limits of sequences.
Thus, $\one\in C$, that is, $X=\one\inv[(1-1/n,1+1/n)]\in\LI(\cU_n)$.

By Proposition \ref{qs}, there are $f_1,f_2,\dots\in A$ such that $\liminf_n f_n\inv[(1-1/n,1+1/n)]=X$.
In particular, $\lim_nf_n=\one$.
\epf

The notation used below is available, e.g., in the survey \cite{ict}.

\bprp
The minimal cardinality of a set $X\sbst\R$ such that $X$ does not
satisfy $\binom{\rmL}{\Ga}$ is $\fb$ (the minimal cardinality of a subset of $\NN$ which is not bounded, with respect to
eventual dominance).
\eprp
\bpf
If $|X|<\fb$, then $X$ satisfies $\sone(\Ga,\Ga)$ \cite{coc2}. Thus, $X$ satisfies $\sone(\Ga,\Ga)$ hereditarily,
and by Lemma \ref{LG}, $X$ satisfies $\binom{\rmL}{\Ga}$.
On the other hand, there is $X\sbst\R$ with $|X|=\fb$, such that $X$ does not satisfy $\sone(\Ga,\Ga)$ \cite{coc2}.
By Proposition \ref{qs}, this $X$ does not satisfy $\binom{\rmL}{\Ga}$.
\epf

The proof of the main theorem in \cite{LinSAdd}, with trivial modifications, gives the first item of the
following theorem. The other items are easy consequences.

\bthm\label{yy}
\mbox{}
\be
\itm
For each unbounded tower $T$ of cardinality $\fb$ in $\roth$, $T\cup\Fin$ satisfies $\binom{\rmL}{\Ga}$.
\itm If $\ft=\fb$, then there are subsets of $\R$ of cardinality $\fb$,
satisfying $\binom{\rmL}{\Ga}$.
\itm There are subsets of $\R$ of cardinality $\ft$, satisfying $\binom{\rmL}{\Ga}$.\qed
\ee
\ethm

The assumption $\ft=\fb$ is known to be strictly weaker than \CH{} or even \MA{}, but it is open
whether it is weaker than $\fp=\fb$, which implies that the sets mentioned in Theorem \ref{yy}
actually have the stronger property $\binom{\Om}{\Ga}$ \cite{LinSAdd}.

\ed